\documentclass[11pt,leqno]{article}
\usepackage[all]{xy}
\usepackage{amssymb,amsmath, amsthm}
\theoremstyle{definition}

\theoremstyle{definition}

\theoremstyle{remark}

\usepackage{amssymb, palatino}
\usepackage{graphicx}

\begin{document}
\title{Theorems of Barth-Lefschetz type in Sasakian geometry}
\author{
 Xiaoyang Chen }
\date{}
\maketitle
\begin{abstract}
{\small In this paper, we obtain theorems of Barth-Lefschetz type in Sasakian geometry.
As corollaries, this implies connectedness principle and Frankel's type theorem.}
\end{abstract}

\noindent {\bf Key words:} {Sasakian manifold \  invariant submanifold \   connectivity}
\\
\noindent{\bf MSC 2010:} { 53C25 53B25 \\   }

\section*{1. Introduction}

\par In the 1920's Lefschetz proved the following theorem now known as the Lefschetz theorem.
Let $N \subset P^k$ be a connected complex submanifold of complex dimension $n$. Let $H$ be a hyperplane and $N\cap H$
a nonsingular hyperplane section. Then the relative cohomology groups satisfy:
\\$$ H^j (N, N\cap H, \mathbb{C})=0, j \leq n-1$$
\par Fifty years later Lefschtz's theorem was generalized by Barth to any two complex submanifolds in $P^k$. More precisely, he proved:
Let $N^n, M^m $ be two connected complex submanifolds of $P^k$, then
$$ H^j(N,N\cap M,\mathbb{C})=0, j \leq min(n+m-k, 2m-k+1).$$
\par Later many people generalized Barth-Lefschetz's theorem,
see [BL],
[FL], [La], [KW], [S1], [S2], [O], [SW].
In [SW], the authors obtained theorems of Barth-Lefschetz type on homotopy groups for complex submanifolds
of certain K\"{a}hler manifolds of nonnegative holomorphic bisectional curvature. As a corollary, this gave a different and elegant proof of
 a well known theorem due to T.Frankel: Suppose $M^{2n}$ is a complete K\"{a}hler manifold with positive
 bisectional curvature and
 $V^{2r}, W^{2s}$ are two compact complex submanifolds with $dimV + dimW \geq dimM $, then $V\cap W \neq\emptyset$.
\par Recently this kind of theorems were
established in other situations. See [FMR], [NW], [Wi] for related work. In this paper, we consider theorems of Barth-Lefschetz type for compact invariant submanifolds of Sasakian manifolds with
nonnegative transversal bisectional curvature. Recall a Riemannian manifold $(M,g)$ is said to be Sasakian if it has a unit Killing field $\xi$,
satisfying the equation
$$ R(X,\xi)Y=<\xi,Y>X-<X,Y>\xi$$
Given such a characteristic vector $\xi$ (also called Reeb vector field), we define a $(1,1)$ tensor $\phi$ to be given by $\phi(X)=\nabla_X \xi $ and the characteristic one form $\eta$
to be given by $\eta(X)=<X,\xi>.$ Altogether we call $(g, \xi, \eta, \phi)$ a Sasakian structure.
The vector field $\xi$ defines the characteristic foliation $F_{\xi}$ with one dimension leaves and the kernel of $\eta$ denoted by $D$, called the contact bundle,
inherits an almost complex structure by restriction of $\phi$.
Let $g^T=g-\eta \otimes \eta$.
It turns out $(g, \xi, \eta, \phi)$ a Sasakian structure iff $(D, g^T, \phi|_D, d\eta)$ defines a transversal K$\ddot{a}$hler structure with transversal K$\ddot{a}$hler form $d\eta$.
\par There is a natural transversal Levi-Civita connection $\nabla^T$ on $M$ defined by
$\nabla_X^T Y=[\nabla_X Y]^p$ if $X,Y\in D$
and $\nabla_\xi^T Y=[\xi, Y]^p$, $\nabla_X ^T \xi=0$, where we denote $Z^p$ is the projection of $Z$ to $D$ for any $Z\in TM$ and
$\nabla$ is the usual Levi-Civita connection induced by $g$.
\par Now define the transversal curvature tensor $R^T$ by
$$ R^T(X,Y)Z=\nabla_X^T\nabla_Y^TZ-\nabla_Y^T\nabla_X^TZ-\nabla^T_{[X,Y]}Z.$$
$$ R^T(X,Y,Z,W)=g<R^T(X,Y)Z,W>.$$
and the transversal bisectional curvature is defined by
$$ B^T(X,Y)=R^T(X,Y,Y,X)+R^T(X,\phi(Y),\phi(Y),X)$$
\\where $X,Y, Z, W \in D$.
\par We say $(M,g)$ has positive (nonnegative) transversal bisectional curvature if and only if $B^T(X,Y)> (\geq) 0$ for any linearly independent vectors $X, Y \in D$
 and $V$ is called an invariant submanifold of $M$ if for any $X\in TV$, $\phi(X) \in TV$.
\par There has been much interest in the classification of compact simply connected Sasakian manifolds with positive (nonnegative)transversal bisectional curvature. See [C], [H], [FOW], [SWZ] for related work. However, as far as the author knows, this program has $not$ been carried out.
\par  Let $M^{2n+1}$ be a Sasakian manifold with contact bundle $D$. Fix $x \in M$ and $X \in D\cap T_xM$.
 Define
 $$ L(x, X)=\{Y \in
  D \cap T_xM: R^T(X,Y,Y,X)+R^T(X,\phi(Y),\phi(Y),X)>0\}$$
  Note $L(x, X)$ is a complex cone and if $Y \in L(x, X)$, then $\phi(Y) \in L(x, X)$. Let $\Gamma$ be the set of complex
 subspaces of $L(x, X)$. Set $l(x, X)= max_{L\in \Gamma} dim_\mathbb{C} (L)$. Define the transversal complex positivity of $M^{2n+1}$
  by $l=inf_{x\in M} l(x)$, where $l(x)=inf_{X \in D\cap T_xM}l(x, X)$. Observe if $M^{2n+1}$ is a Sasakian manifold
  with positive transversal bisectional curvature, then its transversal complex positivity is $n$.
\par The main aim of this paper is to prove the following theorem.
\\
\textbf{ Main Theorem }
{\itshape
Suppose $(M^{2n+1}, g, \xi, \eta, \phi)$ is a complete Sasakian manifold with nonnegative transversal bisectional curvature and $V^{2r+1}, W^{2s+1}$ are two compact invariant submanifolds. If the transversal complex positivity of $M^{2n+1}$ is $l$, then the homomorphism induced by inclusion
$i_*: \pi_j(V, V\cap W)\rightarrow \pi_j(M, W)$ is an isomorhsim for $j\leq r+s-n-(n-l)$ and is a surjection for $j=r+s-n-(n-l)+1$.}
\\
\par We point out our main theorem is $not$ a direct consequence of the results in [SW]! This is explained in section 2.
\par As a corollary of main theorem, we obtain the Sasakian analogue of B. Wilking's connectedness principle, see [Wi].
\\
\\
\textbf{ Corollary 1.1}
{\itshape
Suppose $(M^{2n+1}, g, \xi, \eta, \phi)$ is a complete Sasakian manifold with nonnegative transversal bisectional curvature and has transversal complex positivity
 $l$, then
 \\(1) If $V^{2r+1}$ is a compact invariant submanifold,
the inclusion map $i: V\rightarrow M$ is $(2r-n-(n-l)+1)-$connected;
 \\(2) If $V^{2r+1}, W^{2s+1}$ are two compact invariant submanifolds, $r\leq s$,
the inclusion map $i: V\cap W \rightarrow V$ is $(r+s-n-(n-l))-$connected.
\\Here we say a map $i: V\rightarrow M$ is $k$-connected if the induced map $i_*: \pi_i(V)\rightarrow \pi_i(M)$ is an isomorphism for
$i < k$ and a surjection for $i=k$.}
\\
\par To see part (1), apply main theorem to $V=W$. Now part (2) follows from main theorem and part (1).
\par In section $2$ we will see there are infinite many Sasakian structures with positive transversal bisectional curvature on the standard
sphere $S^{2n+1}$. These are called weighted Sasakian spheres and they contain an interesting class of invariant submanifolds.
They are the so called  Brieskorn manifolds which are codimension two submanifolds of weighted Sasakian spheres. See[B] or [TaTa].
Applying part (1) of corollary 1.1 to  weighted Sasakian spheres, in particular, we obtain the following corollary which is of independent interest.
\\
\\
\textbf{ Corollary 1.2}
{\itshape
 (1) Suppose $V^{2r+1}$ is a compact invariant submanifold of weighted Sasakian spheres. If $r \geq {1 \over 2} (n+1)$, then $V$ is simply connected.
\\
(2)Brieskorn manifolds of dimension $(2n-1)$  are $(n-2)$ connected.}
\\
\par Note the second part of corollary 1.2 was already proved by Milnor, see [M2].
\par As an application of part (2) of corollary 1.1, we get a Frankel's type theorem.
\\
\\
\textbf{ Corollary 1.3}
{\itshape
Suppose $(M^{2n+1}, g, \xi, \eta, \phi)$ is a complete  Sasakian manifold with nonnegative
transversal bisectional curvature and has transversal complex positivity $l$. If
 $V^{2r+1}, W^{2s+1}$ are compact invariant submanifolds with $dimV + dimW \geq dimM +1+2(n-l)$, then $V\cap W \neq\emptyset$ }.
\\
\\In fact, since $dimV + dimW \geq dimM +1+2(n-l)$, then $(r+s)\geq n+(n-l)$. By part (2) of corollary 1.1, the inclusion map $i: V\cap W \rightarrow V$ is $(r+s-n-(n-l))-$connected,
which clearly implies  $V\cap W \neq\emptyset$.
\\
\par Note this Frankel's type theorem was also obtained by other people, see [P1], [TB].
\\
\par The proof of our main Theorem is based on Morse theory on path space.
The crucial part is to choose "good" variational vector fields to get index estimate of critical points of energy functional.
To achieve this, We employ two natural structures on Sasakian manifolds, namely, symmetry structure and transversal K$\ddot{a}$hler structure.
 For details, see section 3.
 \\
 \\
\noindent
{\bf Acknowledgement:} The author would like to express his gratitude to Professor Karsten Grove, who is his advisor and Professor Jianguo Cao,
Professor Xiaoli Cao for helpful advice. He is also grateful to Professor Charles Boyer, Zhi L$\ddot{u}$ and Guofang Wang for discussing examples 2.3
and
2.4. He is also grateful to Weimin Peng for her insistent encouragement.
\section*{2. Examples of Sasakian manifolds and their invariant submanifolds}
In this section we discuss examples of Sasakian manifolds and their invariant submanifolds.
First recall some basic properties of Sasakian manifolds.
\\
\noindent
\textbf{ Proposition 2.1} {\itshape
   Let $(M,g)$ be a Sasakian manifold with Sasakian structure $(\xi,\eta,\phi),$ and $X,Y$ be a pair of
   vector fields on $M$, then
   $$
   \phi^2(X)=-X + \eta(X) \xi $$
   $$\phi(\xi)=0,    <\phi(X), \phi(Y)>= <X,Y> - \eta(X) \eta(Y)
   $$
   $$ (\nabla_X\phi) (Y)= -<X,Y> \xi + <Y, \xi>X $$
   }
 \vspace{5mm}
\noindent
  \textbf {Proposition 2.2 (due to Chinea)} {\itshape
\\ Any invariant submanifold $N$ of Sasakian manifolds is minimal. In fact, let $B$ be
  the second fundamental form of $N$, then $$B(\phi(X), \phi(Y))+B(X,Y)=0.$$
  }
   We refer the reader see [B] for the proof of the above two propositions.
\\
\par A Sasakian structure on $M$ is called quasi regular if the leaves of characteristic foliation $F_{\xi}$ are closed, otherwise it is called irregular. By a theorem of Wadsley [Wa], if a Sasakian structure on $M$ is quasi regular, $\xi$
generates a locally free $S^1$ action on $M$.
A Sasakian structure on $M$ is called regular if this action is free.
If the Sasakian structure on $M$ is quasi regular, then the quotient space $M/F_{\xi}$ is a K$\ddot{a}$hler orbifold. In gerenal, there is no quotient space if it is irregular.
\\
\\
   \textbf{ Example 2.3 (Standard Sasakian sphere) see [BG2]
 }\\
 Let $S^{2n+1}=\{(z_1,\cdots,z_{n+1})\in C^{n+1}|\sum_{k=1}^{n+1}|z_k|^2=1\} $, where
 $z_k=x_k+iy_k ,~x_k,~y_k\in R, ~k=1,2,\cdots,n+1.$
 Equip $S^{2n+1}$ with the Standard metric $g_0$ of constant curvature 1
 and define the Reeb vector field $\xi_0$ by
 $$\xi_0=\sum_{i=1}^{n+1}H_i , $$
 where $H_i=y_i\frac{\partial}{\partial x_i}-x_i\frac{\partial}{\partial y_i}$,
 Then the one form $\eta_0$ takes the form
 $$\eta_0=-\frac{i}{2}\sum_{j=1}^{n+1}(z_jd\overline{z_j}-\overline{z_j}dz_j)$$
 and
 $$\phi_0=\sum_{i,j}((x_ix_j-\delta_{ij})\frac{\partial}{\partial x_i}\otimes dy_j-(y_iy_j-\delta_{ij})\frac{\partial}{\partial y_i}\otimes dx_j
 +x_jy_i\frac{\partial}{\partial y_i}\otimes dy_j-x_iy_j\frac{\partial}{\partial x_i}\otimes dx_j)$$
 \par It is easy to check $S_0=(\xi_0, \eta_0, \phi_0, g_0)$  is a Sasakian structure on $S^{2n+1}$ which we call the standard Sasakian sphere. Also note
 $S^{2k+1}=S^{2n+1}\cap \{(z_1,\cdots , z_{k+1}, 0,\cdots,0)|k\leq n\}$ are invariant submanifolds of the standard sphere.
 \par The standard Sasakian structure on $S^{2n+1}$ is regular. In fact, it is the Hopf fibration as principal $S^1$ bundle over $CP^n$.
 \par \par Using the solution of usual Frankel conjecture, it follows that a simply connected compact Sasakian manifold with
positive transversal bisectional curvature is diffeomorphic to the standard sphere $S^{2n+1}$. However, the next example shows there are infinite many irregular Sasakian structures on
$S^{2n+1}$ with positive transversal bisectional curvature.
 \vspace{5mm}\\
   \textbf{ Example 2.4 (weighted Sasakian sphere) see [BG2] also [TaTa].
 }
 \\Let $w=(w_1, \cdots, w_{n+1})\in R^{n+1}, w_i>0, i=1,2,\cdots, n+1$.
 Define the weighted Sasakian structure $S_w=(\xi_w, \eta_w, \phi_w, g_w)$ on $S^{2n+1}$ by
 $$\xi_w=\sum_{i=1}^{n+1}w_i H_i ,~~~~ \eta_w=\frac{\eta_0}{\sum_{i=1}^{n+1}w_i((x_i)^2+(y_i)^2)}$$
 $$\phi_w=\phi_0-\frac{1}{1+\eta_0(\rho_w)}\phi_0(\rho _w)\otimes \eta_0$$
 $$g_w=d\eta_w\circ(\phi_w\otimes1)+\eta_w\otimes\eta_w , $$
 where $\rho_w=\xi_w-\xi_0. $
 \par T.Takahashi first showed $S_w$  are Sasakian structures on $S^{2n+1}$. They are deformations of the standard Sasakian sphere which
  are called weighted Sasakian spheres. Observe if $w=(w_1, \cdots, w_{n+1})$ is a $(n+1)$-tuple of rational numbers, there is a natural orbifold fibration $S^{2n+1}\rightarrow CP(w)$,
 where $CP(w)$ is the weighted projective space. Hence in this case, the weighted Sasakian structures on $S^{2n+1}$ are quasi regular and have positive transversal bisectional
 curvature since weighted projective spaces are global quotients of $CP^n$. However, if $w=(w_1, \cdots, w_{n+1})$  are linearly independent over rational numbers,
 the corresponding Sasakian structures are irregular. In this case, one can approximate it by a family of quasi regular
 weighted Sasakian structures. Hence we obtain infinite many irregular Sasakian structures on $S^{2n+1}$ with positive transversal
 bisectional curvature! These examples show
our main theorem is $ \mathbf{not}$ a direct consequence of results in [SW]!
By a conversation with Professor G. Wang, the author knows it is conjectured any compact simply connected Sasakian manifold with positive transversal
 bisectional curvature is transversally biholomorphic to a weighted Sasakian sphere. See [BG2] for the notion of transversally biholomorphic.
 \par Now we discribe a very interesting class of invariant submanifolds of weighted Sasakian spheres, see[B] or [TaTa].
 Let $w=(w_1, \cdots, w_{n+1})$ be a $(n+1)$-tuple of positive integers.
 A polynomial of the form
 $ P(z)=z_1^{w_1}+\cdots+z_{n+1}^{w_n+1}, (z_1, \cdots, z_{n+1})\in C^{n+1}$  is called  a Brieskorn polynomial. Then
  $ \sum^{2n-1}(w_1,\cdots,w_{n+1})=P^{-1}(0)\cap S^{2n+1}(1)$ is called a Brieskorn manifold. It is a codimension two submanifold
  of $S^{2n+1}$. T. Takahashi first showed it is an invariant submanifold of weighted Sasakian sphere [TaTa].
 \\
 \\  \textbf{ Example 2.5 ( see [ToTa])}
 Let $G$ be a simply connected simple Lie group and $H$ be a compact Lie subgroup. Suppose $G/H$ is an irreducible Hermitian symmetry space.
 Assume $g, h$ are the Lie algebra of $G$ and $H$, respectively. Now let
 $$g=h+p$$
  be the canonical decomposition with respect to $B$, where $B$ is the Killing form of $g$.
  Let $(J, Q)$ be the K\"{a}hler structure on $G/H$ such that $Q_0=\frac{-1}{8n}B_p$ if $G$ is compact and $Q_0=\frac{1}{8n}B_p$
  if G is noncompact. Here $Q_0$ is the restriction of $Q$ to the tangent space $T_0(G/H)$ of $G/H$ at $0=[H]\in G/H$ which is identified
  with $p$, $B_p$ is the restriction of $B$ to  $p$ and $2n$ is the dimension of $G/H$. It is known that the center $z(h)$ of $h$ is one dimension and there exists $Z_0\in z(h)$
  such that
  $$ J_0=-ad_pZ_0,$$
where $J_0$ is the restriction of $J$ to the tangent space $T_0(G/H)$, $ad_p$ is
the restriction to $p$ of the adjoint representation of $g$ in $g$ ([KN]).
Since $H$ is compact, we have the direct sum decomposition
$$h=[h,h]+z(h)$$
Now let $H_0$ be a connected Lie subgroup of $H$ with Lie algebra $[h,h]$. By a theorem of Kato-Motomiya ([KM]),
$G/H_0$ admits a $G$-invariant Sasakian structure $(g, \xi, \eta, \phi)$
such that
$$\phi=-ad_k Z_0, \xi =2Z_0, \eta _0=\frac{1}{2}Z_0^*,$$
and
$g=-\frac{1}{8n}B_k$ if $G$ is compact or
$g=\frac{1}{8n}B_k + 2\eta_0\otimes \eta_0$, if $G$ is noncompact.
\\
Here $k=p+z(h)$ is identified with the tangent space $T_0(G/H_0)$,
$Z_0^*$ is a $1$-form on $k$ defined by $Z_0^*(Z_0)=1$ and $Z_0^*(X)=0$ for all $X\in p$.
\par The Sasakian manifold $G/H_0$ is a principal circle bundle over Hermitian symmetry space $G/H$.
It is a Sasakian $\phi$-symmetry space which is defined in [ToTa] and is a regular Sasakian manifold.
Hence if $G$ is compact,  $G/H_0$ has nonnegative transversal bisectional curvature
since $G/H$ is a Hermitian symmetry space of compact type. Moreover,
the transversal complex positivity of $G/H_0$ is equal to the complex positivity of $G/H$ and the latter is computed in
[KW].
\par Let $\pi:G/H_0\rightarrow G/H$ be the projection map, then for any complex submanifold $N$ of $G/H$,
$\pi^{-1}(N)$ is an invariant submanifold of $G/H_0$.
 \par In K\"{a}hler geometry, due to a uniformation theorem of N. Mok ([M3]), up to taking covers, any compact K\"{a}hler manifold with
 nonnegative bisectional curvature is a product of flat space and hermitian symmetry space. However, due to the knowledge of the author,
 he doesn't know any classification result or conjecture on compact Sasakian manifolds with nonnegative transversal
 bisectional curvature.

\section*{3. Morse theory and index estimate of critical points}
\noindent
\noindent
\par We follow the same setup as in [SW]. Suppose $M^{2n+1}$ is a complete Sasakian manifold with nonnegative transversal bisectional curvature and $V^{2r+1}, W^{2s+1}$ are two compact invariant submanifolds. Also assume the transversal complex positivity of $M^{2n+1}$  is $l$.
Let $\Omega_{V,W}$ be the piecewise smooth path space of $M$ with starting point in $V$ and endpoint in $W$.
Consider the energy functional on $\Omega_{V,W}$
\\$$ E(\gamma)=\int^1_0|\dot{\gamma}|^2 dt$$
The critical points of $E$ are geodesics that start perpendicularly to $V$ and end perpendicularly to $W$.
Let $\gamma$ be such a geodesic. The index of $\gamma$ is closely related to the topology of $\Omega_{V,W}$.
In fact, by Morse theory, we have the following
\\
\\
\textbf {Lemma 3.1}
{\itshape Suppose $V$ intersects $W$ transversally and every nontrivial critical points of $E$ on
$\Omega_{V,W}$ has index $\lambda> \lambda_0 \geq 0$. Then the relative homotopy groups $\pi_i(\Omega_{V,W}, V\cap W)$ are zero for
$0\leq i \leq \lambda_0$.}
\\
\\
\begin{proof} See [SW] section 1.
\end{proof}
\par Now we apply an idea of B.Wilking and P. Petersen to choose "good" variational vector fields along $\gamma$. The idea is to employ the symmetry structure of
$M$ in nature. More precisely, choose vector fields $X_1, X_2, \cdots, X_r$ along $\gamma$ such that
$$ X_i(0)\in T_{\gamma(0)}V $$
$$ \dot{X_i}= - <X_i, \nabla_{\dot{\gamma}}\xi>\xi $$
and $X_1(0), X_2(0), \cdots, X_r(0), \phi (X_1)(0), \phi(X_2)(0), \cdots, \phi(X_r)(0), \xi(0)$
are orthogonal basis of $T_{\gamma(0)}V$.
By taking derivative along $\gamma$, it is easy to see
$$ <\xi, \dot{\gamma}>=0$$
$$ <\xi, X_i>=0$$
$$ <X_i,\dot{\gamma}>=0$$
along $\gamma$.
Furthermore, since $V$ is an invariant submanifold, it is easy to see
$\phi(\dot{\gamma}(0)) \bot{T_{\gamma(0)}V}$ and $\phi(\dot{\gamma}(1)) \bot{T_{\gamma(1)}W}$.
Although $X_i$ are not parallel vector fields along $\gamma$, it is a surprising fact that $\phi(X_i)$ are, where $i$ always ranges from 1 to $r$
in this section.
\\
\vspace{5mm}
   \noindent
   \textbf{Lemma 3.2.} {

  $$ \nabla_{\dot{\gamma}}(\phi(X_i))=0 $$
  $$ <\dot{\gamma}, (\phi(X_i))>=0$$
    }

\vspace{5mm}

\begin{proof}
$$ \nabla_{\dot{\gamma}}(\phi(X_i))=(\nabla_{\dot{\gamma}}\phi)(X_i))+
\phi ( \nabla_{\dot{\gamma}}(X_i))=$$
$$-<\dot{\gamma}, X_i> \xi + <X_i,\xi> \dot{\gamma}+ \phi(-<X_i,\nabla_{\dot{\gamma}}\xi>\xi)=0,
$$
where the second equality follows from the facts $\nabla_{\dot{\gamma}}<\dot{\gamma}, \phi(X_i)>=0 $
and $<\dot{\gamma}, \phi(X_i)>(0)=0$.
\end{proof}
\textbf{Remark 3.3:} It is easy to see $\nabla^T_{\dot{\gamma}} X=0$. Hence $X$ is parallel with respect to the transversal Levi-Civita connection.
\par Now we perturb the metric using Cheeger's deformation. Denote $G$ be the one parameter subgroup generated by the Killing vector
field $\xi$ and consider the diagonal action of $G$ on $M\times G $. For arbitrary positive number $\lambda$, the product metric $ g+\lambda dt^2 $
on $M\times G$
induces a metric $g_{\lambda}$ on $M$ such that $Q: M\times G\longrightarrow M $ defined by $(p,x)\rightarrow x^{-1}p$ is a Riemannian submersion.
Under this new metric on $M$, $\gamma$ is still a geodesic and the length of $\xi$ tends to zero as $\lambda$ goes to zero. However, directions orthogonal to $\xi$ remain unchanged and the curvatures $R(\dot{\gamma},X_i,X_i, \dot{\gamma}), R(\dot{\gamma},\phi(X_i),\phi(X_i),\dot{\gamma})$
can only become larger as $X_i, \phi (X_i), \dot{\gamma}$ are perpendicular to $\xi$. Furthermore, by direct computation it is not hard to see $\|\bigtriangledown_{\dot{\gamma}}(X_i)\|_{g_\lambda}\rightarrow 0$ as $\lambda\rightarrow 0$, see [P2] page 100.
\par Now denote $\nabla$, $\nabla^{\lambda}$ the Levi-Civita connection induced by $g$, $g_{\lambda}$ and $R$, $R^{\lambda}$ the usual Riemannian curvature tensor with respect to $g$ and $g_{\lambda}$, respectively. Let $R^T$ be the
transversal Riemannian curvature tensor with respect to the Sasakian metric $g$, then we have the following
\\
\textbf{ Lemma 3.4}
$$ R^T(X,Y,Y,X)=R(X,Y,Y,X)+3|g<X,\phi(Y)>|^2, where X,Y\perp \xi.$$
\begin{proof}
$$ R^T(X,Y,Y,X)=g<R^T(X,Y)Y,X>=g<\nabla_X^T\nabla_Y^TY-\nabla_Y^T\nabla_X^TY-\nabla^T_{[X,Y]}Y,X>$$
$$ =g<\nabla_X(\nabla_YY-<\nabla_YY,\xi>\xi)-\nabla_Y(\nabla_XY-<\nabla_XY,\xi>\xi)$$
$$-\nabla^T_{[X,Y]-<[X,Y],\xi>\xi}Y-\nabla^T_{<[X,Y],\xi>\xi}Y,X>$$
$$=R(X,Y,Y,X)+<\nabla_XY,\xi><\nabla_Y\xi,X>+<[X,Y],\xi><\nabla_Y\xi,X>$$
$$=R(X,Y,Y,X)+3|g(X,\nabla_Y\xi)|^2,$$
\\where the last equality follows from $X, Y \perp \xi$ and $\xi$ is a Killing vector field.
\end{proof}
\par In view of lemma 3.2 and lemma 3.4, it is easy to see
\\
\textbf{Lemma 3.5}
Suppose $(M, g)$ is a Sasakian manifold with positive (nonnegative) transversal bisectional curvature, then
$$R(\dot{\gamma}, X_i, X_i, \dot{\gamma})+ R(\dot{\gamma}, \phi(X_i), \phi(X_i), \dot{\gamma})>0 (\geq0)$$
\\
$$R^{\lambda}(\dot{\gamma}, X_i, X_i, \dot{\gamma})+ R^{\lambda}(\dot{\gamma}, \phi(X_i), \phi(X_i), \dot{\gamma})>0 (\geq0)$$
\\
\textbf{ Lemma 3.6}
 $$ <\nabla_X^{\lambda}X,\dot{\gamma}>_{g_\lambda}|_0^1+<\nabla_{\phi(X)}^{\lambda}\phi(X),\dot{\gamma}>_{g_{\lambda}}|_0^1=0,
 where  X\perp \xi. $$
\begin{proof}
Since
$$\phi(X)\bot\xi$$
$$X\bot\xi,$$
$$ \dot{\gamma} \bot\xi$$
  $(X,0), (\phi(X), 0), (\dot{\gamma}, 0)$ are horizontal vector fields on $(M\times G, g+\lambda dt^2)$.
 Denote $\overline{X}$ be the horizontal lift of $X$. Using the fact $ \nabla_{\overline{X}}{\overline{Y}}={\overline{\nabla_X^{\lambda}Y}}+ \frac{1}{2}[{\overline{X}},{\overline{Y}}]^v$ (see [P2] page 94), we see
$$\nabla_{\overline{X}}{\overline{X}}={\overline{\nabla_X^{\lambda} X}}$$
and so $$ <\nabla_X^{\lambda}X,\dot{\gamma}>_{g_\lambda}|_0^1+<\nabla_{\phi(X)}^{\lambda}\phi(X),\dot{\gamma}>_{g_{\lambda}}|_0^1$$
$$=<{\overline{\nabla_X^{\lambda} X}},\overline{\dot{\gamma}}>_{g+\lambda dt^2}|_0^1+<\overline{\nabla_{\phi (X)}^{\lambda} \phi (X)},\overline{\dot{\gamma}}>_{g+\lambda dt^2}|_0^1$$
$$ =<\nabla_{\overline{X}}{\overline{X}},\overline{\dot{\gamma}}>_{g+\lambda dt^2}|_0^1+<\nabla_{\overline{\phi(X)}}\overline{\phi(X)},\overline{\dot{\gamma}}>_{g+\lambda dt^2}|_0^1$$
$$= <\nabla_X X,\dot{\gamma}>_g|_0^1+<\nabla_{\phi(X)}\phi(X),\dot{\gamma}>_g|_0^1=0$$
where the last equality follows from proposition 2.2.
\end{proof}
\par Now we compute the second variation of $X_i$ and $\phi (X_i)$ with respect to $g_\lambda$.
$$ {1\over 2} E_{**}(X_i,X_i)= <\nabla_{X_i} X_i, \dot{\gamma}>|_0^1+ \int_0 ^1< \nabla_{\dot{\gamma}}(X_i), \nabla_{\dot{\gamma}}(X_i)>-R(\dot{\gamma}, X_i, X_i, \dot{\gamma}) dt$$
$$ {1\over 2} E_{**}(\phi(X_i), \phi(X_i))= <\nabla_{\phi(X_i)} \phi(X_i), \dot{\gamma}>|_0^1+ \int_0 ^1< \nabla_{\dot{\gamma}}(\phi(X_i)), \nabla_{\dot{\gamma}}(\phi(X_i))>$$
$$-R(\dot{\gamma}, \phi(X_i), \phi(X_i), \dot{\gamma}) dt$$
From the above discussions,  we get
$${1\over 2} E_{**}(X_i,X_i)+{1\over 2} E_{**}(\phi(X_i), \phi(X_i))=<\nabla_{X_i} X_i, \dot{\gamma}>|_0^1+ <\nabla_{\phi(X_i)} \phi(X_i), \dot{\gamma}>|_0^1
$$
$$ +\int_0 ^1< \nabla_{\dot{\gamma}}(X_i), \nabla_{\dot{\gamma}}(X_i)>+< \nabla_{\dot{\gamma}}(\phi(X_i)), \nabla_{\dot{\gamma}}(\phi(X_i))>$$
$$-R(\dot{\gamma}, X_i, X_i, \dot{\gamma}) -R(\dot{\gamma}, \phi(X_i), \phi(X_i), \dot{\gamma}) dt$$
$$=\int_0 ^1< \nabla_{\dot{\gamma}}(X_i), \nabla_{\dot{\gamma}}(X_i)>
-R(\dot{\gamma}, X_i, X_i, \dot{\gamma})-R(\dot{\gamma}, \phi(X_i), \phi(X_i), \dot{\gamma}) dt$$

$$\rightarrow \int_0 ^1 -R(\dot{\gamma}, X_i, X_i, \dot{\gamma}) -R(\dot{\gamma}, \phi(X_i), \phi(X_i), \dot{\gamma}) dt \leq 0$$
as $\lambda \rightarrow 0$.
\par To get strict inequality we want
$$\int_0 ^1 R(\dot{\gamma}, X_i, X_i, \dot{\gamma}) +R(\dot{\gamma}, \phi(X_i), \phi(X_i), \dot{\gamma}) dt > 0$$
This is insured by requiring that
$$X_i, \phi(X_i) \in L \bigcap \xi^\bot_{\gamma (0)} $$
where $L$ is a complex subspace of $L(\gamma(0), \dot{\gamma}(0))$ such that
$dim_\mathbb{C} L=l(\gamma(0), \dot{\gamma}(0))$
and $\xi^\bot_{\gamma (0)}$ is the orthogonal complement of $\xi$ in $T_{\gamma(0)}V$. Let
$P=L\bigcap \xi^\bot_{\gamma (0)}$. Then we have
$dim_\mathbb{C}P\geq l+r-n$.
Parallel transport $P$ along $\gamma$ to $\gamma(1)$ and denote the resulting subspace by $K$.
Then $K \bigcap \xi^\bot_{\gamma (1)}$ has complex dimension at least $r+s-n-(n-l)+1$,
where  $\xi^\bot_{\gamma (1)}$ is the orthogonal complement of $\xi$ in $T_{\gamma(1)}W$.
So for sufficient small $\lambda$, we can find a metric $g_\lambda$ on $M$ such that the index of $\gamma $ is at least
$\geq r+s-n-(n-l)+1$. So we have shown
\\
\\
\textbf{Theorem 3.7}
{\itshape Suppose $(M^{2n+1}, g, \xi, \eta, \phi)$ is a complete Sasakian manifold with nonnegative
transversal bisectional curvature and $V^{2r+1}, W^{2s+1}$ are compact invariant submanifolds.
If the transversal complex positivity of $M^{2n+1}$ is $l$, then
for every nontrivial critical point $\gamma$ of the energy functional $E$ on $\Omega_{V,W}$, we have
$$index (\gamma)\geq r+s-n-(n-l)+1$$
 and so by lemma 3.1
 $$\pi_i(\Omega_{V,W}, V\cap W)=0,
0\leq i \leq r+s-n-(n-l)$$. }
\section*{4. Proof of main theorem}
\noindent
\noindent
 Fix $x\in V$,  let $\Omega^W_x$ be the space of piecewise smooth pathes that start from $x$ and end in $W$. Consider the fibration:
 \[
\xymatrix{
  \Omega^W_x \ar[r]^{}
                & \Omega_{V,W} \ar[d]^{e}  \\
                & V             }
\]
where $e$ takes $\gamma$ to $\gamma (0)$.
It is well known $ \pi_i(\Omega^W_x) = \pi_{i+1}(M,W)$.
From the long exact sequence of homotopy groups for fibration, we have
$$ \rightarrow\pi_{i+1}(\Omega_{V,W})\rightarrow \pi_{i+1}(V)\rightarrow \pi_{i+1}(M,W)\rightarrow \pi_i (\Omega_{V,W})\rightarrow \pi_i(V)\rightarrow$$
is exact.
\\
Consider the diagram
\[
\xymatrix{
   \pi_{i+1}(V\cap W)\ar[d]_{} \ar[r]^{}
                & \pi_{i+1}(V) \ar[d]_{\simeq} \ar[r]^{}& \pi_{i+1} (V, V\cap W) \ar[d]^{} \ar[r]^{} & \pi_i(V\cap W)\ar[d]_{} \ar[r]^{}
                & \pi_i(V) \ar[d]_{\simeq}   \\
  \pi_{i+1}(\Omega_{V,W}) \ar[r]_{}
                & \pi_{i+1}(V) \ar[r]_{} & \pi_{i+1}(M,W) \ar[r]_{} & \pi_i(\Omega_{V,W}) \ar[r]_{} & \pi_i(V)
                  }
\]
where the vertical arrows are induced by inclusion. So
\\$$i_*: \pi_j(V, V\cap W)\rightarrow \pi_j(M, W)$$ is an isomorphism for $j\leq r+s-n-(n-l)$ and is a surjection for $j=r+s-n-(n-l)+1$ by theorem 3.7
and the commutivity of the diagram.

\bigskip

\noindent
{
Xiaoyang Chen \\
{\small Department of Mathematics}\\
{\small University of Notre Dame} \\
{\small 46556,Notre Dame, Indiana, USA}\\
{\small 574-631-3741} \\
{\small xchen3@nd.edu}

\end{document}